\documentclass{article}

\usepackage{amsmath}

\newtheorem{definition}{Definition}
\newtheorem{lemma}[definition]{Lemma}

\newtheorem{proposition}[definition]{Proposition}

\newcommand{\qed}{\hfill \rule{1ex}{1ex}} 

\newenvironment{proof}{{\bf Proof}: }{\qed}

\title{Restricted permutations and queue jumping}
\author{
M.~H.~Albert\footnotemark[1]
\and
R.~E.~L.~Aldred\footnotemark[2]
\and
M.~D.~Atkinson\footnotemark[1]
\and
H.~P. van Ditmarsch\footnotemark[1]
\and
C.~C.~Handley\footnotemark[1] 
\and
D.~A.~Holton\footnotemark[2]
}

\begin{document}
\footnotetext[1]{Department of Computer Science, University of Otago}

\footnotetext[2]{Department of Mathematics and Statistics, University 
of Otago}
\maketitle

\begin{abstract}
A connection between permutations that avoid $4231$ and a certain queueing discipline is established.  It is proved that a more restrictive queueing discipline corresponds to avoiding both $4231$ and $42513$, and enumeration results for such permutations are given.
\end{abstract}

\section{Introduction}
Let $s_n(\alpha_1,\alpha_2,\ldots)$ be the number of permutations of length $n$ that avoid all of the permutation patterns $\alpha_1,\alpha_2,\ldots$.  Finding a formula or a generating function for $s_n(\alpha_1,\alpha_2,\ldots)$ is a difficult and much studied problem.  For a single pattern $\alpha$ the sequence $s_n(\alpha)$ is known only in the following cases (see \cite{Bona2,Gessel,Knuth}):

\begin{itemize}
\item $\alpha=12\ldots k$ or $\alpha=k\ldots 21$ for any $k$
\item $|\alpha|\leq 3$
\item $|\alpha|=4$ but $\alpha\not=1324, 4231$
\end{itemize}

So the first unsolved cases are $\alpha=1324$ and $\alpha=4231$ which are equivalent by symmetry.  Here the lower bound $s_n(4231)\geq s_n(4321)$ has been proved by B\'ona \cite{Bona1} who has also made some contributions towards an upper bound.

In this note we introduce a certain type of queue and, harking back to early work by Knuth \cite{Knuth} on such problems, study its connection with $4231$-avoiding permutations.  By strengthening the conditions satisfied by the queue we go on to solve the enumeration problem for the sequence $s_n(4231,42513)$.

\section{Jump queues}

We shall define two queue-based data structures. They will be used with an input sequence $1,2,\ldots,n$ whose members are added one by one to the rear of the queue.  Removals from the queue generate an output sequence that will be a permutation of $1,2,\ldots,n$.  Removals are always allowed from the front of the queue but both our queue-based structures also allow elements other than the front member of
the queue to be output, that is {\em queue jumping} (when we refer to
a {\em jump} we will always mean an output operation which would not
be permitted by an ordinary queue). When a queue jump occurs some of the items in the queue become ``locked'' (forbidden to jump until the lock is released).

Our two structures differ in
the extent to which instances of queue jumping restrict further queue
jumping. The first of these data structures, the {\em loosely locked
jump queue} is defined by the property that when an element $x$ is jumped from the queue all the elements in the queue behind $x$ become locked; they are released (given the freedom to jump) when all the elements in front of $x$ have been output.  Note that any new elements added to a loosely locked jump queue are initially not locked.

In the second structure, the {\em strictly locked jump queue}, the locking rule is more severe.  Again, elements in the queue behind an element $x$ that jumps become locked; furthermore, any new elements that are added to a queue that already has some locked elements, are initially in the locked state.  As before, the lock on an element is released once all the elements in front of the jumped item that initially caused the lock have been output.

In both cases jumping from the rear of the queue imposes no locks so our queues are at least as powerful as the input-restricted deques that were analysed in \cite{Knuth}.

In studying the output permutations generated by either type of queue we observe that an output permutation may be producible in several ways.  This allows us to make a simplifying assumption. Suppose that we are
attempting to generate a particular permutation $\pi$ and have
proceeded to a point where we wish to output a symbol $p$. If $p$ is
already in the queue then if it is possible to succeed at all from
this point onwards, we can succeed by outputting $p$ immediately. For
the only alternative is to add further elements to the queue and then
output $p$. The only effect that this might have (versus outputting
$p$ immediately and then adding the same elements) is to lock some
queue elements that would not be locked in the original instance. So, it
cannot be harmful to do any output as soon as it becomes
available and, from now on, we consider only operation sequences with this property.  Under that assumption we will regard the production of any permutation $\pi$ as taking place in a number of \emph{stages}.  In any of these stages one or more input elements are added to the queue, the last of these is then output (such outputs produce the left to right maxima of $\pi$), and then further output from the queue occurs (possibly none at all); a stage comes to an end when the next element of $\pi$ to be output has not yet been added to the queue.

We begin our investigation with a result whose easy proof is omitted.

\begin{lemma}Suppose we have a jump queue of either sort with a frontal segment $\alpha=a_1<a_2<\ldots <a_m$.  Then the permutations of $\alpha$ that can be generated by queue removals (from the front or by jumping) are precisely those that avoid $231$.
\end{lemma}
\begin{proposition}
The collection of permutations that can be produced by a loosely
locked jump queue is the class of $4231$-avoiding permutations.
\end{proposition}

\begin{proof}
Let a permutation $\pi$ be given which contains a $4231$ pattern as 
\[
\pi = \cdots d \cdots b \cdots c \cdots a \cdots
\]
and suppose that we could produce $\pi$ using a loosely locked jump
queue. In order to output $d$ before all of $a$, $b$, and $c$, those
elements must be in the queue when $d$ is output. However, the
subsequent output of $b$ would then lock $c$ so that it could not be output
until $a$ was. So, in fact, $\pi$ could not be generated.

Conversely, suppose that a permutation $\pi$ avoids $4231$. Let $m_1<m_2<\ldots<m_k$ be the left to right maxima of $\pi$.  Then we can write
\[
\pi = m_1 \, \alpha_1 \, m_2 \, \alpha_2 \, \cdots \, m_k \alpha_k
\]
where each $\alpha_i$ is some segment of $\pi$ and $m_i > \alpha_i$. For
convenience define $m_0 = 0$.

We will show that $\pi$ can be produced by following the operation of
the queue in attempting to produce it, and observing that we never
reach a point where an element which we need to output is locked.  We argue inductively on the  stages of this production (as defined previously) where stage $j$ produces the segment $m_j\alpha_{j}$.

In the first stage the elements from $1$ through $m_1$ are added to the
queue,  then $m_1$ is output (without causing any locks), and then the elements of
$\alpha_1$ must now be output. At present this can certainly be
accomplished since $\alpha_1$ avoids $231$. However, this sequence of
operations may leave some remaining elements of $[1, m_1)$ locked in
the queue. 

Suppose that, after $j$ stages,  we have output the initial segment $m_1  \alpha_1
\cdots m_j \alpha_j$ of $\pi$. Stage $j+1$  begins by adding
the elements of $(m_j, m_{j+1}]$ to the queue and then outputting the
element $m_{j+1}$. Next we begin to output the elements of
$\alpha_{j+1}$. Consider the point at which an element $c$ of this
type is to be output. Choose $i \leq j+1$ such that $m_{i-1} < c <
m_i$. Any lock to the output of $c$ would have been applied by an 
element $b < c$ jumping after $c$ had been added to
the queue (i.e. after the output of $m_i$). For this lock to have
remained in force there must be an element $a < b$ still in the
queue. But if all this were true the elements $m_i b c a$ would form a
$4231$ pattern in $\pi$.  Stage $j+1$ therefore succeeds in producing a further segment $m_{j+1}\alpha_{j+1}$ of output and the inductive proof is complete.
\end{proof}

\begin{proposition}\label{sljq}
The collection of permutations that can be produced by a strictly
locked jump queue is the class of $\{4231, 42513\}$-avoiding permutations.
\end{proposition}

\begin{proof}
The proof is similar in spirit to that above, so we will provide
somewhat fewer details. In one direction, suppose that $\pi$ contains the pattern $42513$ as the subsequence $dbeac$ and yet can be produced by the queue.  Then the element $c$ would be locked by the output of $b$.  Since it is required that $a$ be output after $e$ the element $a$ must still be in the queue when $e$ enters and so $c$ would still be locked at this point; therefore $e$ would be locked upon entering the queue and could not be output at the proper place.
Obviously if $\pi$
contains $4231$ it cannot be output, for even a loosely locked queue
would not suffice in that case.

Conversely, suppose that $\pi$ avoids these two patterns. Write
\[
\pi = m_1 \, \alpha_1 \, m_2 \, \alpha_2 \, \cdots \, m_k \alpha_k
\]
as above, and follow the operation of the queue in stages again.

Suppose that the output of $m_j$ is prevented because it has been locked by some preceding jump of
an element $b$. Then $b$ was preceded by some $m_i > b$, and some element $c >
b$ still remains in the queue as does some element $a < b$ (or $m_j$
would no longer be locked). Then one of the two sequences $m_i b m_j c
a$ or $m_i b m_j a c$ occurs in $\pi$. The latter is a $42513$ pattern and
the former contains a $4231$ pattern, a contradiction in either case.

Now suppose that the output of some $x$ in $\alpha_j$ is prevented by a lock caused by
some preceding jump of an element $b$. Then choose $m_i$, $c$, and $a$
as above. Since the pattern $4231$ does not occur in $\pi$ the elements $m_i$, $b$, $x$, $a$, $c$
occur in $\pi$ in that order. If $x > m_i$ we have a pattern $42513$ while if $x < m_i$
then $m_i b x a$ is a $4231$ pattern. Again, a contradiction is achieved in
either case.
\end{proof}

\section{The number of permutations avoiding $4231$ and $42513$}

Consider the operation of a strictly locked jump queue as it produces some permutation that avoids $4231$ and $42513$. At a point
where no elements are locked we might choose to add one or more input elements (after
which our next output step must be to remove the last element of the
queue), to jump an element from the rear of the queue (which imposes
no locks), or to output an earlier element, say the $j$th. In the
latter case we must, or rather might as well, output all elements (if any)
which are earlier still in the queue before continuing the operation. In so
doing, we can produce any $231$-avoiding permutation of these $j-1$
elements. As is well known, the number of such permutations is
$c_{j-1}$, the $(j-1)$st Catalan number.

If we set $q$ to be the number of elements in the queue, and $i$ the
number remaining in the input then this trichotomy is easily
translated into a recurrence for the number of permutations that can be produced from this configuration. However, the manipulation of the
resultant quantities will be simplified if we make a distinction
between two cases: where the next output is the last element of the
queue, or where we do not place any restriction on the next
output. Let $l(q,i)$ enumerate the former class and $n(q,i)$ the
latter (both quantities conventionally $0$ if either $q$ or $i$ is
negative). Then we have:
\begin{eqnarray*}
l(q,i) &=& l(q+1,i-1) + n(q-1,i) \\
n(q,i) &=& l(q,i) + \sum_{j=1}^{q-1} c_{j-1} n(q-j, i).
\end{eqnarray*}

Consider the power series:
\begin{eqnarray*}
N(x,t) &=& \sum_{q, i \geq 0} n(q,i) x^q t^i \\
L(x,t) &=& \sum_{q, i \geq 0} l(q,i) x^q t^i \\
C(x) &=& \sum_{i \geq 0} c_i x^i \\
f(t) &=& L(0,t),
\end{eqnarray*}
the last of which, by Proposition \ref{sljq}, is the generating function for the class
of permutations that avoid $4231$ and $42513$. Then the recurrences translate easily into the
equations:
\begin{eqnarray*}
L(x,t) &=& \frac{t(L(x,t) - f(t))}{x} + xN(x,t) + 1 \\
N(x,t) &=& L(x,t) + x N(x,t) C(x) - x f(t) C(x).
\end{eqnarray*}
After some rearrangement this yields:
\[
L(x,t) \left( (x^2 - x t)C(x) + x^2 - x + t \right) =
\left( t - x t C(x) + x^3 C(x) \right) f(t) -x + x^2 C(x).
\]
Now the stage is set for a simple application of the kernel
method. Consider the circumstances under which the
parenthetical expression in the left hand side above is zero:
\begin{eqnarray*}
 (x^2 - x t)C(x) + x^2 - x + t &=& 0 \\
x C^2(x) - C(x) + 1 &=& 0.
\end{eqnarray*}
Solving formally for $C(x)$ and $x$ in terms of $t$ yields:
\begin{eqnarray*}
C^3(x) t - C(x) + 1 &=& 0 \\
x &=& C(x) t
\end{eqnarray*}
These conditions can then be substituted in the right hand side, which
must also be zero, yielding eventually:
\[
f(t) = \frac{1}{1 - \eta(t) t}
\]
where
\[
t \eta(t)^3 - \eta(t) + 1 = 0.
\]
The latter equation is easily seen to be the equation satisfied by the generating function for ternary trees.  Therefore the $n$th term of the power series expansion of $\eta(t)$ is the number of ternary trees on $n$ nodes, namely
\[\frac{\binom{3n}{n-1}}{n}\]
From this it follows readily that, if $f_n$ denotes the $n$th term of $f(t)$, 
\[\liminf_{n\rightarrow\infty} \sqrt[n]{f_n}=27/4\] 
and more detailed asymptotics could easily be obtained.  From the above equations we can also read off the recurrence
\[f_n=\sum_{i=0}^{n-1}f_{n-1-i}\binom{3i}{i-1}/i\]
and thereby compute the expansion of
\[f(t)=1+t+2t^2+6t^3+23t^4+102t^5+495t^6+2549t^7+13682t^8+75714t^9+\ldots\]
to as many terms as necessary.
 
\end{document}